# Robust adaptive strategies for the guidance of users in road networks

Farida Manseur[1], Nadir Farhi, Habib Haj-Salem, Jean-Patrick Lebacque

*Université Paris-Est, IFSTTAR/COSYS/GRETTIA, F-77447 Champs-sur Marne Cedex France.*

**Abstract**

We present an algorithm for optimal guidance of users in road networks. It is a "stochastic-on-time-arrival (SOTA)"-like algorithm which calculates optimal guidance strategies with reliable paths, for road network origin-destination pairs. Our contribution consists here in extending an existing SOTA algorithm, in order to include robustness of the guidance strategy, towards path failures.

The idea of SOTA algorithms is to calculate the maximum probability of reaching a destination node, starting from any node of a road network, and given a time budget. This calculus gives the optimal path for every origin-destination pair of nodes in the network, with an associated optimal adaptive guidance strategy, and with respect to the considered time budget. The approach models very well the importance of the travel time variability in the route choice process. Indeed, by adopting this kind of guidance strategies, users may accept short deteriorations on the average travel time, if in return, they have a guarantee on the reliability of the travel time (in term of the maximum probability of reaching their destinations in given time budgets).

We propose here an extension of this approach in order to take into account the existence and the performance of alternative detours of the selected paths, in the calculus of the guidance strategy. We take into account the fact that one or many links of the selected optimal path may fail during the travel. We then consider that users may be sensitive to path changing. That is to say that they may prefer paths with efficient alternative detours, with respect to paths without, or with less efficient detours, even with a loss in the average travel time, and/or in its reliability. In order to take into account such behaviors, we propose a model that includes the existence as well as the performance of detours for selected paths, in the calculus of the travel time reliability (i.e. the maximum probability of reaching a destination node). This new way of calculating travel time reliability guarantees a kind of robustness of the guidance strategies. That is to say that, the travel time reliability associated to the obtained optimal guidance strategy is not likely to change, however associated adaptive paths change during the travel. The variation of the travel time reliability, with respect to a network structure changing, is thus improved.

We show the effectiveness of our algorithm on a grid network, where the vector of the travel times through the links of the network is assumed to have a multivariate gamma probability distribution, taking into account the correlations of the travel times. We present preliminary results of the application of the proposed guidance strategy, with comparisons to the existing algorithm. We show, in particular, the variance of the calculated optimal paths and strategies, with respect to the users' sensitivity towards the existence and the performance of alternative detours of the selected paths. The model is then able to calculate robust guidance strategies towards path failures.

## 1. Introduction

The optimal guidance problem of users in transportation networks knows a renewed interest in recent years. This is due mainly to the fast development of new technologies of information and communications. This development has made possible and necessary the development of real-time mobile and/or GPS computing applications, where adaptive algorithms give optimal solutions, updated according to the traffic conditions, and with better guarantees in terms of reliability, robustness and risk, compared to prior solutions. Optimal routing in transportation networks with highly varying traffic conditions is a challenging problem due to the stochastic nature of travel-time on links of the network. Depending on the users' preferences, an optimal route might need

---

*Corresponding author. E-mail address:* farida.manseur@ifsttar.fr

to consider notions of both the expected value the travel-time and its reliability (variance). This is a multi criterion optimization problem that is in general hard to solve.

Loui (1983) defined a formulation of routing optimality as a route with the Least Expected Time (LET). The LET problem has been well investigated and numerous efficient algorithms have been proposed to solve the different variants of the problem; see for example Fu and Rilett (1998), Miller-Hooks and Mahmassani (2001), waller and Ziliaskopoulos (2002). However, there are several frameworks in which the LET solution is not adequate, because it does not take into consideration the variance of travel-time distributions and does not give any guarantee on its reliability.

A very natural definition of a reliable optimal path is given in the Stochastic On Time Arrival (SOTA) approach, presented by Frank (1969). This formulation permits yet to obtain satisfactory adaptive solution in terms of travel time reliability. A formulation of the SOTA problem using stochastic optimal control is presented by Bertsekas (2005). Fan and al (2006) formulated it as a stochastic dynamic programming problem, and solved it using a standard successive approximation (SA). However, in presence of cycles on the network, as is the case with all road networks, there is no finite bound on the maximum number of iterations required for the convergence of the algorithm. As an alternative, Nie et al. (2006) proposed a discrete approximation algorithm for the SOTA problem which converges in a finite number of steps and runs in pseudo-polynomial time.

Samaranayake et al. (2012) presented a number of optimization techniques in order to speed up the computation time of the algorithm. The approach includes a label-setting algorithm based on the existence of a uniform strictly positive minimum link travel-time, advanced convolution methods centered on the Fast Fourier Transform and the approch of zero-delay convolution, and localization techniques for determining an optimal order of policy computation. However, the performance of all of these algorithms is limited by the large search space of the problem. Sabran et al. (2014) have shown how preprocessing methods can be used to further reduce the computation time of the SOTA problem. Unfortunately, the structure of the SOTA problem formulation limits the types of preprocessing methods that can be used for this problem, and prevents massive running time reductions in the deterministic case. Recently, Kobitzsch et al. (2014) presented a novel approach to reduce the immense computational effort of stochastic routing based on existing techniques for alternative routes. In an extensive experimental study, they showed that the process of stochastic route planning can be speed-up immensely without sacrificing much in terms of accuracy.

In this paper, we focus on the reliable path problem that aims to find an adaptive optimal routing strategy, which takes into account the existence and the performance of alternative detours of the selected paths in the road network. That is to say, we assume that one or many links of the selected optimal path may fail during the travel, and that users may be sensitive to path changing. Some users may then prefer paths with efficient alternative detours, with respect to paths without or with less efficient detours.

Our approach is based on the idea of Frank (1969) considering that a reliable path from a given origin to a given destination maximizes the probability of realizing a travel time less than a given time budget. Thus, we propose an adaptation of this approach introducing robustness in the selection of the optimal path. We base here on the routing model of Samaranayak (2011). From the probability distributions of travel times through the links of the network, users evaluate their maximum probability to reach their destination in given time budgets, and through different possible routes. The model takes into account correlations between the travel times through the links of the network.

The remainder of this paper is organized as follows. In section 2, the SOTA problem and the SA method are summarized. The formulation of the problem and the description of our proposed algorithm are provided in section 3. Computational tests are conducted in section 4. In section 5 conclusions are drawn.

## 2. Summary of the SOTA problem formulation and SA algorithm

In this section, we briefly summarize the SOTA problem for the convenience of readers and the continuity of our discussion in later sections. The details of the approach summarized here are available for example in Samaranayak et al. (2011).

We consider the graph $G(N, A)$, where $N$ is the set of nodes, with $|N| = n$, and $A$ is the set of arcs, with $|A| = m$. The set of successor and predecessor nodes are denoted by $\Gamma^{+1}(i) = \{j|(i,j) \in A\}$ and $\Gamma^{-1}(i) = \{k|(k,i) \in A\}$ respectively. The SOTA problem is to find the best routing strategy from any starting node $i$ ($i = 1, 2, \ldots, d$) maximizing the probability of arriving at the destination node $d$ within a desired time $T$. The routing procedure from a given origin to a given destination is considered as a multistage decision process, every decision stage of the routing process includes two states; the physical state which represents the location of the current decision point, and the information state which represents the time budget for completing the remaining

process. The optimal decision for the remaining process is made on the basis of the two states of the current stage. This represents a typical feedback control process in which dynamic programming plays a primary role. Due to the randomness involved in the transition between steps, and contrary to deterministic processes, optimal strategy can not be identified a priori.

Given a node $i \in N$ and a time budget $t$, $u_{ki}(t, y)$ denotes the maximum probability of reaching destination node $d$ parting from origin node $i$, within time $t$, under the optimal policy, where $k$ is the predecessor node from which the user is arriving, and $y$ is the realized travel-time on the upstream link $(k, i)$. At each node $i$, the traveler should choose the link $(i, j)$ that maximizes the probability of arriving on time at the destination. According to Samaranayake and al (2011), the maximum probabilities $u_i(t, k, y)$ satisfy:

$$u_{ki}(t, y) = \max_{j \in \Gamma^{+1}} \int_0^t p_{ij}(tt_{ij} = w \mid tt_{ki} = y) u_{ij}(t - w, w) dw \tag{1}$$

$$\forall i \in N \setminus \{d\}, k \in \Gamma^{-1}(i), j \in \Gamma^{+1}(i), 0 \leq t \leq T, 0 \leq y \leq T - t$$

$$u_{kd}(t, y) = 1, \quad \forall k \in \Gamma^{-1}(d), 0 \leq t \leq T, 0 \leq y \leq T - t, \tag{2}$$

where $tt_{ij}$ denotes the travel time on link $(i, j)$, and $p_{ij}(tt_{ij} \mid tt_{ki})$ is the probability distribution function (pdf) of $tt_{ij}$, conditioned on $tt_{ki}$. The pdf $p_{ij}(.)$ is assumed to be known and can for example be obtained using historical data or real-time traffic information. $u_{ij}(t - w, w)$ is the maximum probability of arriving to destination node $d$ within time $t - w$, parting from node $j$, conditioned that the travel time on link $(i, j)$ is $w$.

One approach to solve this problem is to use a *successive approximation* (SA) method (**Algorithm 1**) proposed by Fan and Nie (2006). Algorithm 1 solves the system equation (1) iteratively until convergence, and gives an optimal routing policy

Given initial approximations, the iterative relationships for successive approximation are given by Algorithm1.

| **Algorithm 1**. Successive approximation algorithm (Fan and Nie (2006)) |
|---|
| **Step 0**. Initialization<br>$iter = 0$ (iteration index)<br>$u_{ki}^{iter}(t, y) = 0, \quad \forall i \in N \setminus \{d\}, k \in \Gamma^{-1}(i), 0 \leq t \leq T, 0 \leq y \leq T - t,$<br>$u_{kd}^{iter}(t, y) = 1, \quad \forall k \in \Gamma^{-1}(d), 0 \leq t \leq T, 0 \leq y \leq T - t,$<br><br>**Step 1**. Update<br>$iter = iter + 1$<br>$u_{kd}^{iter}(t, y) = 1, \quad \forall k \in \Gamma^{-1}(d), 0 \leq t \leq T, 0 \leq y \leq T - t,$<br>$u_{ki}^{iter}(t, y) = \max_{j \in \Gamma^{+1}} \int_0^t p_{ij}(tt_{ij} = w \mid tt_{ki} = y) u_{ij}^{iter-1}(t - w, w) dw$<br>$\hspace{3cm} \forall i \in N \setminus \{d\}, k \in \Gamma^{-1}(i), j \in \Gamma^{+1}(i), 0 \leq t \leq T, 0 \leq y \leq T - t$<br><br>**Step 2**. Convergence test<br>If $\forall (i, t) \in N \times [0, T], \max_{i,t} |u_{ki}^{iter}(t, y) - u_{ki}^{iter-1}(t, y)| = 0,$ then stop;<br>Otherwise go to step 1. |

At each iteration $iter$, $u_{ki}^{iter}(t, y)$ gives the probability of reaching the destination node $d$ from origin node $i$ within a time budget $t \in [0, T]$, using a path that contains no more than $iter$ links, and under the optimal policy. A convergence property is proposed by Fan and Nie (2006) using the bounded monotone convergence theorem. According to this theorem, the approximation error monotonically decreases with $iter$ and the solution reaches an optimal value when $iter$ is equal to the number of links in the optimal path. However, optimal paths in general stochastic networks may include loops. In this case, the number of iterations required to converge is undetermined.

## 3. Robust guidance

In this article we base on the routing model (equations (1) and (2)) presented by Samaranayake and al. (2011), where from the probability distributions of travel times through the links of the network, users evaluate their maximum probability to reach their destination in given time budget, and using different possible routes. The model takes into account correlations between the travel times through the links.

We propose here an extension of this approach in order to take into account the existence and the performance of alternative detours of the selected paths, in the calculus of the guidance strategy. We take into account the fact that one or many links of the selected optimal path may fail during the travel. We then consider that users may be sensitive to path changing. That is to say that they may prefer paths with efficient alternative detours, with respect to paths without, or with less efficient detours, even with a loss in the average travel time, and/or in its reliability. In order to take into account such behaviors, we propose a model that includes the existence as well as the performance of detours for selected paths, in the calculus of the travel time reliability (i.e. the maximum probability of reaching a destination node). This new way of calculating travel time reliability guarantees a kind of robustness of the guidance strategies. That is to say that, the travel time reliability associated to the obtained optimal guidance strategy is not likely to change, however associated adaptive paths change during the travel. The variation of the travel time reliability, with respect to a network structure changing, is thus improved. For that, we propose to calculate for each node $i$ the maximum probability $u_{ki}(t, y)$ to reach the destination node, where we take into account the case where the selected path fails before the users who selected it reach the destination node; for witch case, alternative neighboring paths are used. $u_{ki}(t, y)$ denotes, as above, the maximum probability to reach the destination node, departing from node $i$, in a time budget $t$, and knowing that the user comes from node $k$ upstream of $i$, and that the realized travel time from $k$ to $i$ is $y$. However, the mathematical definition of $u_{ki}(t, y)$ is different from (1)-(2), and is given below.

3.1. *Calculation of the probabilities $u_{ki}(t, y)$*

We consider the following notations:

$$A_{kij}(t, y) = \int_0^t p_{kij}(tt_{ij} = w \mid tt_{ki} = y) u_{ij}(t - w, y) dw$$

$$\forall i \in N \backslash \{d\}, k \in \Gamma^{-1}(i), j \in \Gamma^{+1}(i),$$

$$\forall\ 0 \leq t \leq T, 0 \leq y \leq T - t.$$

We then sort, for a given $(k, i)$ pair, the quantities $A_{kij}(t, y)$ on $j$, in a decreasing order, and denote them $B_{kip}(t, y)$ in that order, i.e. $B_{ki1}(t, y) \geq B_{ki2}(t, y) \geq \cdots \geq B_{ki\lambda}(t, y)$. We then have

$$\mathbb{A}_{ki}(t, y) := \{A_{kij}(t, y), \forall j \in \Gamma^{+1}(i)\} = \{B_{kip}(t, y) \in \mathbb{A}_{ki}(t, y), B_{ki1}(t, y) \geq B_{ki2}(t, y) \geq \cdots \geq B_{ki\lambda}(t, y)\}$$

Then we write the maximum probability for a user to reach destination node from node $i$ in a time budget $t$, knowing that the user comes from node $k$ and that the travel time from $k$ to $i$ is $y$, as follows.

$$u_{ki}(t, y) = \sum_{p=1}^{\lambda} \Psi_p B_{kip}(t, y), \quad \forall\ i \neq d, k \in \Gamma^{-1}(i), 0 \leq t \leq T, 0 \leq y \leq T - t \tag{3}$$

$$u_{kd}(t, y) = 1, \quad \forall\ k \in \Gamma^{-1}(d), 0 \leq t \leq T, 0 \leq y \leq T - t \tag{4}$$

with $\Psi_p$ are weighting coefficients satisfying.

$$\sum_{p=1}^{\lambda} \Psi_p = 1$$

Where, in case where $\lambda > |\Gamma^{+1}(i)|$, $B_{kip}(t,y) = 0$ for $p > |\Gamma^{+1}(i)|$. In this case, $\lambda$ is the number of successors to be taken account in the network, independent of $i$. Indeed, one can distinguishes the following two cases:
- Case 1 : $\lambda < |\Gamma^{+1}(i)|$, in wich case, only some successors of $i$ are considered in the mean.
- Case 2 : $\lambda > |\Gamma^{+1}(i)|$, in wich case, we have $\sum_{j \in \Gamma^{+1}} \Psi_j(t) < \sum_{p=1}^{\lambda} \Psi_p(t) = 1$.

In the case 2 above, nodes $i$ that have a small number of successors are penalized; they get low values $u_{ki}(t,y)$. Therefore, paths that pass through these nodes i.e. paths with small number of alternatives or detours shall have low probabilities to be selected as optimal paths. One way to choose $\lambda$ would be to take the maximum over the cardinals of the sets $\Gamma^{+1}(i)$ of successors of all the nodes of the network.

$$\lambda = \max_{i \in N} |\Gamma^{+1}(i)|,$$

where $|\cdot|$ denotes the cardinal of a set.

In order that formula (3) has meaning, $\Psi_p$ have to be chosen such that $\Psi_1 \geq \Psi_2 \geq \cdots \geq \Psi_\lambda$. That is to say that $\Psi_p$ decrease as $B_{kip}(t,y)$ decrease on $p$. This dependence of $\Psi_p$ on $A_{kip}(t,y)$ makes the model nontrivial.

Let us notice here that if $\lambda = 1$ or if $\lambda > 1$ and $\Psi_p = 0, \forall p \geq 2$, then (3)-(4) coincides with (1)-(2).

*3.2. Calculation of successor nodes*

In the calculus of $u_{ki}(t,y)$ with *(3)*, instead of maximizing the quantities $A_{kij}(t,y)$, we propose to take a wiegted mean of theses quantities, with weights $\Psi_p$, $p = 1,2,\ldots\lambda$.

The optimal guidance strategy is then determined by the sequence of successor nodes $s_{ki}(t,y)$ as shown in the following formula.

$$s_{ki}(t,y) = \arg \max_{j \in \Gamma^{+1}(i)} (A_{kij}(t,y)) \quad (5)$$

$s_{ki}(t,y)$ denotes here the optimal successor node of node $i$ for u user to reach the destination node, knowing that the user comes from node $k$ upstream of $i$, and that the realized travel time on link $(k,i)$ is $y$. By taking a mean in (3) rather than the maximum (as in (1)), we do not only take into account the path maximizing the probabilities $u_{ki}(t,y)$, but we also take into account the existence and the performance of alternative deviations at each node. The SA algorithm for the formulation (3)-(5) is summarized as follow:

---
**Algorithm 2**.

**Step 0**. Initialization
$\lambda = \max_{i \in N} |\Gamma^{+1}(i)|$ (maximum number of successor nodes to take into account)
$iter = 0$ (iteration index)
$u_{ki}^{iter}(t,y) = 0, \quad \forall i \in N \backslash \{d\}, k \in \Gamma^{-1}(i), 0 \leq t \leq T, 0 \leq y \leq T-t,$
$u_{kd}^{iter}(t,y) = 1, \quad \forall k \in \Gamma^{-1}(d), 0 \leq t \leq T, 0 \leq y \leq T-t,$

**Step 1**. Update
$iter = iter + 1$
$u_{kd}^{iter}(t,y) \triangleq 1, \quad \forall k \in \Gamma^{-1}(d), 0 \leq t \leq T, 0 \leq y \leq T-t,$
$u_{ki}^{iter}(t,y) = \sum_{p=1}^{\lambda} \Psi_p B_{kip}^{iter}(t,y) \quad \forall i \neq d, k \in \Gamma^{-1}(i), 0 \leq t \leq T, 0 \leq y \leq T-t$
$s_{ki}^{iter}(t,y) = \arg \max_{j \in \Gamma^{+1}(i)} (A_{kij}^{iter}(t,y)) \quad \forall i \neq d, k \in \Gamma^{-1}(i), 0 \leq t \leq T, 0 \leq y \leq T-t$

---

In the following, we study the efficiency of our algorithm (Algorithm 2) by applying it on the grid network of 25 nodes and 40 links, shown in Figure 1.

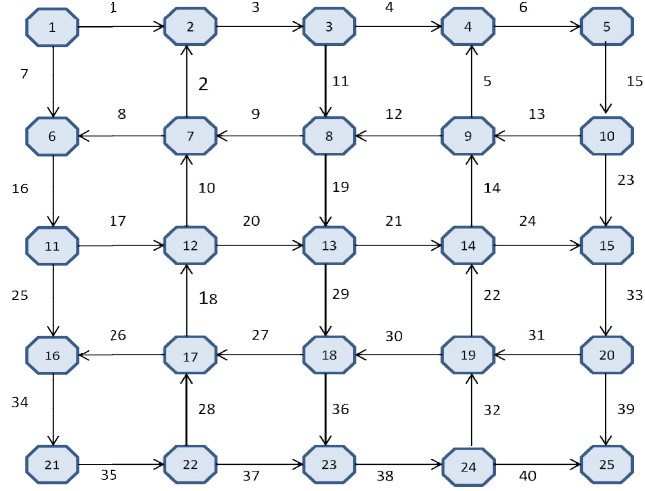

Figure 1. A 25-node network with gamma distributed link travel-time

We assume that link travel times on the network of Figure 1 are drawn from bi-variate Gamma distribution. More precisely, we assume that the joint probability distribution of two successive links of the network is a bi-variate Gamma. We base here on the bi-variate Gamma distribution given by Smith et al. (1982). Indeed, the joint PDF of two positively correlated random variables $W$ and $Y$, with Gamma marginal distributions $f_W$ and $f_Y$ with shape and scale parameters $(\alpha_W, \beta_W)$ and $(\alpha_Y, \beta_Y)$ respectively, was derived by Smith et al. (1982) as follow.

$$f(w,y) = \begin{cases} \dfrac{g_1}{g_2} \sum_{k_1=0}^{\infty} \sum_{k_2=0}^{\infty} c_{k_1 k_2} (\beta_W w)^{k_1} (\eta \beta_Y y)^{k_1+k_2} & \text{if } \rho > 0 \\ f_W(w) \cdot f_Y(y) & \text{if } \rho = 0 \end{cases} \qquad (6)$$

where $\rho$ is the product-moment correlation coefficient of $W$ and $Y$, estimated from the sample data,

$$\rho = \frac{E[(W - m_W)(Y - m_Y)]}{\sigma_W \sigma_Y}$$

$m_W, m_Y$ and $\sigma_W, \sigma_Y$ are the sample means and standard deviations of the variables $W$ and $Y$ respectively.

$$w, y \geq 0, \ 0 < \eta = \rho\sqrt{\alpha_y/\alpha_w} < 1, \ \alpha_y \geq \alpha_w, \ 0 \leq \rho \leq \eta\sqrt{\alpha_w/\alpha_y}$$

$$g_1 = (\beta_W w)^{\alpha_W - 1} (\beta_Y y)^{\alpha_Y - 1} \exp\left(-\frac{\beta_W w + \beta_Y y}{1 - \eta}\right) \qquad (7)$$

$$g_2 = (1 - \eta)^{\alpha_W} \Gamma(\alpha_W) \Gamma(\alpha_Y - \alpha_W) \qquad (8)$$

$$c_{k_1 k_2} = \frac{\eta^{k_1+k_2} \Gamma(\alpha_Y - \alpha_W + k_2)}{(1 - \eta)^{2k_1 + k_2} \Gamma(\alpha_Y + k_1 + k_2) k_1! k_2!} \qquad (9)$$

$\Gamma(.)$ is the Gamma function $\Gamma(z) = \int_0^\infty t^{z-1} e^{-t} dt$.

## 4. Numerical example

In the network of Figure 1, the maximum number of successors over all the nodes is equal to 2. Therefore, according to Section 3, we take here $\lambda = 2$. Then, we have two weighting coefficients $\Psi_1$ and $\Psi_2$. We simply denote here $\Psi = \Psi_1$, and then $\Psi_2 = 1 - \Psi$.

To reach destination node 25 departing from node 1, we have 154 elementary paths. We apply Algorithm 2 and derive the probabilities $u_{ki}(t, y)$ for all origin nodes $i$ of the network. In order to illustrate our approach, let us consider the probability $u_{23}(t, y)$ of reaching the destination node 25 from node 3. We chose origin node 3 here because from node 3 we have two routing actions: go to successor node 4, which will give us only one routing action at the next step (go to successor node 5), or go to successor node 8, which will give us two routing actions at the next step (go to 7 or 13). Indeed, if we chose node 4 as successor of node 3, then at node 4, if link 6 fails, the whole routing strategy fails with it. However, if we chose node 8 as successor of node 3, then at node 8, if one of the links 9 or 19 fails, we have the possibility to change and take a detour using the other link; see Figure 1.

Then we consider the following scenario. Travel-times on two successor links of the network follow a bi-variate gamma probability distribution, with given average travel time vector and variance-covariance matrix. We take the same average travel time (9 time units) for all the links, except link 4, for which we take a travel time of 5 time units. We take, for every pair of two successive links, the same covariance matrix (variance = 3, covariance = 1.5). We assume that we are at node 3 and that we have spent a time $y = 9$ time units on link 3, upstream of node 3. The maximum time budget we consider here is 80 time units.

Figure 2 shows the maximum probabilities $u_{23}(t, 9)$ in function of time budget $t$, and for different values of $\Psi$. Figure 3 gives the optimal successor nodes $s_{23}(t, 9)$ in function of time budget $t$, and for different values of $\Psi$.

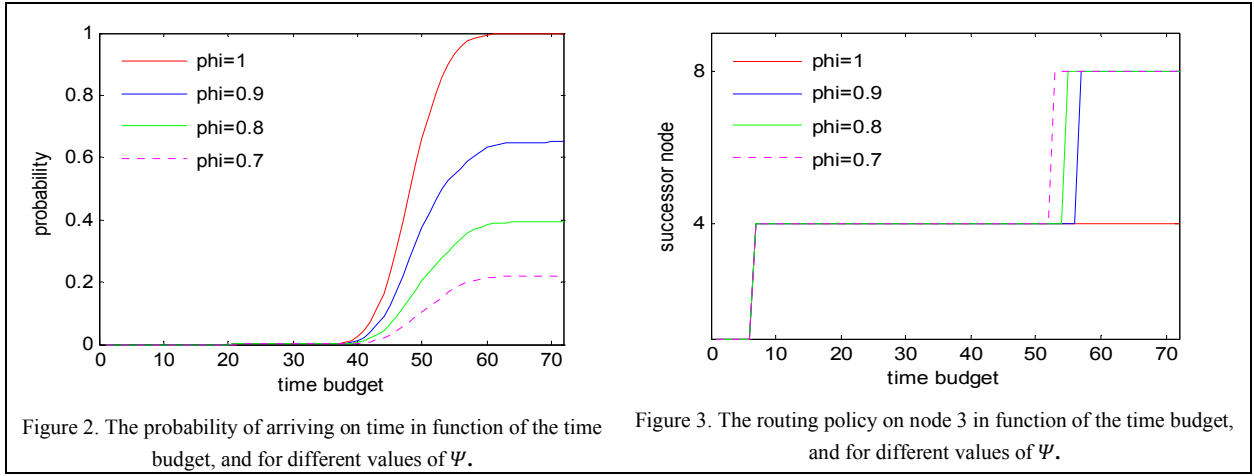

Figure 2. The probability of arriving on time in function of the time budget, and for different values of $\Psi$.

Figure 3. The routing policy on node 3 in function of the time budget, and for different values of $\Psi$.

The case $\boldsymbol{\Psi = 1}$ corresponds to the model (1)-(2), while the cases $\boldsymbol{\Psi < 1}$ correspond to the model (3)-(4). We see from Figure 2 that the maximum probability $u_{ki}(t, y)$ decreases as the values of $\boldsymbol{\Psi}$ decrease. That is to say that for lower values of $\boldsymbol{\Psi}$, lower maximum probabilities $u_{ki}(t, y)$ are obtained. This is because we replaced a maximum operator in (1)-(2) by a mean value in (3)-(4). Indeed, a user taking a lower value of $\boldsymbol{\Psi}$, asks for more path-robustness or path-flexibility, and, in the counterpart, he loses in term of travel time reliability. The difference $u_{ki}(t, y, \boldsymbol{\Psi_1}) - u_{ki}(t, y, \boldsymbol{\Psi_2})$ can then be interpreted as the price of path-robustness corresponding to a measure of it given by the difference $\boldsymbol{\Psi_1 - \Psi_2}$.

Figure 3 gives the optimal routing policy in term of the optimal successor node. We see from that figure, that for $\Psi = 1$, the optimal successor node is node 4, for all the time budgets up to 80 time units. Indeed, two successor nodes (4 and 8) are possible for node 3. However, since the average travel time on link 4 (from node 3 to node 4) is (5 time units) lower than the average travel time (9 time units) on link 11 (from node 3 to node 8), and since all the other parameters are symmetric in the network (average travel times, variances, etc.), and since we do not take into account path-robustness here ($\Psi = 1$), it is trivial that paths passing from node 4 are better than those passing by node 8. With low values of $\Psi = 0.9$ (respectively 0.8, 0.7), where path-robustness is considered, we see that, with a time budget bigger than 57 (respectively 55, 53) time units, the optimal successor node of node 3 is node 8 rather than node 4, even though paths passing through node 4 have lower average travel time comparing to those passing through node 8. For $\Psi = 0.9$ for example, with a time budget of 57 time units, paths passing through node 8, with travel time reliability of $u_{23}(57,9) = 0.59$ are preferable than those passing from node 4, with a travel time reliability of $u_{23}(57,9) = 0.97$; see red and blue curves in Figure 2 and Figure 3.

For $\Psi = 0.9$, in general, we see that with a time budget lower than 56 time units, the algorithm prefers paths passing through node *4* (see blue line in Figure 3). However, with a time budget higher than 56 time units, the optimal policy on node *3* changes, and node *8* becomes the optimal successor node (see blue line in Figure 3). That means that, node 4 which has only one successor is penalized, i.e. it gets low values $u_{ki}(t, y)$. Therefore, paths that pass through that node i.e. paths with small number of alternatives or detours have low probability to be selected as optimal paths.

From Figure 2, if we want to reach the destination node within a time budget of *60* time units, then we get a travel time reliability of *0.99* with $\Psi = 1$, and a travel time reliability of *0.63* with $\Psi = 0.9$**.** Then, the value $0.36 = 0.99 - 0.63$ can be interpreted as the *price of robustness* to pay in term of travel time reliability. On the other hand, if we want to reach the destination node with a travel time reliability of *0.60*, we need a time budget of *49* time units with $\Psi = 1$ and a time budget of *58* time units with $\Psi = 0.9$. Then, the value $9 = 58 - 49$ can be interpreted as the *price of robustness* to pay in term of travel time budget.

From these results we can conclude that:
- If a user maximizes the travel time reliability of the paths, without taken into account their robustness, then he will choose the path that passes from node 4, because it is the one maximizing the probability of reaching the destination node 25 in the considered time budget.
- If a user seeks a guarantee in terms of robustness even if he loses in terms of travel time budget, then he will choose the path that passes from node 8.

## 5. Conclusion

This article considers the routing problem of maximizing the probability of arriving on time with a new robust adaptive strategy for the guidance of users in road networks. An existing model is modified. In order to include the performance of alternative detours of the selected paths, we extended the concept of reliability by introducing a new reliability index. The modification that we made to the existing model allows the selection of an optimal path according to two criteria: the reliability of the path in term of travel time and the robustness of the path in term of flexibility (existence and performance of alternative detours).